\documentclass[11pt]{amsart}

\usepackage{amsthm}
\usepackage{amssymb}
\usepackage{amsmath}

\newcommand{\labell}[1] {\label{#1}}

\numberwithin{equation}{section}

\newtheorem {Theorem}   {Theorem} 
\numberwithin{Theorem}{section}

\newtheorem {Lemma}[Theorem]    {Lemma}         
\newtheorem {Proposition}[Theorem]{Proposition}  
\theoremstyle{definition}
\newtheorem{Definition}[Theorem]{Definition}
\theoremstyle{remark}
\newtheorem{Remark}[Theorem]{Remark}

\newtheorem {Corollary}[Theorem]{Corollary}  

\newtheorem {Claim}[Theorem]{Claim}

\def	\N {{\mathbb N}}
\def    \haf {\frac{1}{2}}

\def	\C {{\mathbb C}}

\def	\R	{{\mathbb R}}

\def  \ra {\rightarrow}

\begin{document}


\setlength{\smallskipamount}{6pt}
\setlength{\medskipamount}{10pt}
\setlength{\bigskipamount}{16pt}





\title[Periodic Orbits of Hamiltonian Flows Near Symplectic Minima]{Periodic Orbits of Hamiltonian Flows Near Symplectic Critical Submanifolds}

\author[Ely Kerman]{Ely Kerman}

\address{Department of Mathematics, UC Santa Cruz, 
Santa Cruz, CA 95064}
\email{ely@cats.ucsc.edu}

\date{October, 1998}

\thanks{This work was partially supported by the NSF}

\bigskip

\begin{abstract}

In this paper we produce a lower bound for the number of periodic orbits of certain Hamiltonian vector fields near Bott-nondegenerate symplectic critical submanifolds. This result is then related to the problem of finding closed orbits of the motion of a charged low energy particle on a Riemannian manifold under the influence of a magnetic field.
\end{abstract}

\maketitle

\section{Introduction} 
Our primary goal is to establish a lower bound for the number of periodic orbits of a low energy, charged particle moving in a magnetic field on a Riemannian manifold $M$. This motion can be described by a Hamiltonian dynamical system on $T^*M$ with the metric Hamiltonian and the ``twisted'' symplectic form $\Omega=d\lambda+\pi^* \omega.$ Here $\lambda$ is the standard Liouville 1-form, $\pi \colon T^*M \ra M$ is the canonical projection and $\omega$ is a closed 2-form on $M$ representing the magnetic field. We will show that in the low energy (large field) limit, when $\omega$ is nondegenerate, the flow restricts to the fibres of $T^*M$ and is quasiperiodic there. Hence, one might expect many of the periodic orbits to persist on low energy levels. For a large class of cases, we prove that this happens and establish a lower bound for the number of these orbits. To the knowledge of the author, it is still unknown whether this number can indeed be finite.

When $M$ is a surface, the problem of finding periodic orbits has been studied extensively. For a comprehensive review of these results the reader is referred to \cite{gi:Cambr}. Symplectic techniques were first applied to the problem on the torus by V.I. Arnold in \cite{Arn1} and \cite{Arn2}. These methods were refined by V.L. Ginzburg and dividing the search for periodic orbits into the three natural realms of high, intermediate and low energy levels we can summarize his results as follows. On high energy levels, when the magnetic flow is close to the geodesic flow, he established the existence of a periodic orbit in each free homotopy class of maps $S^1 \ra M$ that contains a closed geodesic, \cite{gi:Glik}. The only nontrivial case here, when $M$ is the torus, is dealt with using a theorem due to J. Mather, \cite{ma}. In \cite{gi:MathZ}, Ginzburg noted that the horocycle flow of Hedlund \cite{hed} provides an example of a magnetic flow on an intermediate energy level with no closed trajectories. Most importantly for us, Ginzburg showed that on low energy levels the number of closed orbits is no less than $Crit(M)$, the minimum number of critical points of a smooth real-valued function on the surface $M$, \cite{gi:FA}, \cite{gi:MathZ}.

Much less is known about the existence of periodic orbits when $dim(M)>2$. In \cite{gi:ker} it is shown that for any magnetic field on a torus of any dimension there exist periodic orbits on almost all energy levels. This is proven using the symplectic capacity of Hofer and Zehnder \cite{hz} and the results concerning this capacity found in \cite{fhv} and \cite{jiang}.  For nonzero magnetic fields on tori, L. Polterovich \cite{pol} showed, using Hofer's metric on the space of Hamiltonian diffeomorphisms \cite{hof1}, that there actually exist contractible closed orbits on a sequence of arbitrarily small energy levels. As well, A. Bahri and I. A. Taimanov have used variational methods to prove a specialized existence result for exact magnetic fields on manifolds of any dimension, \cite{bata}.

When $M$ is even dimensional and $\omega$ is nondegenerate as well as closed, $M$ is a Bott-nondegenerate symplectic critical submanifold of $T^*M$. Thus, the search for low energy periodic orbits can be placed in the more general setting of looking for periodic orbits of a Hamiltonian dynamical system near such submanifolds. This is the setting in which we will work. Letting $d^2_N H$ denote the restriction of the Hessian of $H$ to vectors transversal to $M$, and requiring the eigenvalues of $d^2_N H$ to satisfy certain resonance conditions with respect to $\Omega$, we establish the existence of at least $qCL(M)+(n-m)$ closed orbits on sufficiently low energy levels. Here $q$ is the number of integer independent eigenvalue classes of $d^2_N H$ with respect to $\Omega$ and $CL$ denotes the cuplength. When $M$ is a point the lower bound becomes $n$, in agreement with A. Weinstein \cite{we:sta}, and our method of proof reduces to that of J. Moser \cite{Moser-1976}. For the magnetic problem, we also get that if the metric Hamiltonian comes from an almost comlex structure compatible with $\omega$, then $q=1$ and there are at least $CL(M)+m$ closed orbits on low energy levels. This includes the Kahlerian metric on any Kahler manifold. It is plausible that the dependence of our lower bound on $q$ is an artifact of the particular method we employ to detect periodic trajectories and the real bound should be $(n-m)(CL(M)+1).$ 

In looking for periodic orbits, the limiting dynamics leads us to consider the following situation. Let $X_0$ be a vector field whose flow contains a periodic compact nondegenerate submanifold $\Sigma$. For us $\Sigma$ will be nondegenerate if the linearized Poincar\'{e} map has 1 as an eigenvalue of algebraic multiplicity $dim(\Sigma)$. Given a perturbation $X_1$ of $X_0$, we would like to know under what conditions $X_1$ has periodic orbits near $\Sigma$. Moreover, we would like to obtain a lower bound for the number of these closed orbits. To answer these questions we consider M. Bottkol's generalizations \cite{Bottkol} of the related work of Moser \cite{moser:kep} \cite{Moser-1976} and Weinstein \cite{we:persub} \cite{we:sta}. In this work, Bottkol  produces a lower bound for the number of periodic orbits of $X_1$ near $\Sigma$, given that both $X_0$ and $X_1$ are Hamiltonian with respect to the same symplectic form.

Bottkol's method of proof is extended to a setting, which covers our limiting dynamics, in which $X_0$ is not Hamiltonian. With this, the stated lower bound for the number of closed orbits of certain Hamiltonian flows near Bott-nondegenerate symplectic critical submanifolds is established. This result is discussed in Section \ref{sec:theorem} as is its relation to the magnetic problem. The proof, as well as a tangential result concerning closed characteristics of certain 2-forms, is contained in Section \ref{sec:proof}. 

\subsection*{Acknowledgments.}I am profoundly grateful to Viktor Ginzburg for showing me this problem, helping in the solution and his kind thoughtful advice. I would also like to thank Alan Weinstein for his helpful comments.

\section{Statement and Discussion of the Main Theorem}
\labell{sec:theorem}

\subsection{The Main Theorem}

Let $H \colon N^{2n} \ra {\R}$ be a proper function on a symplectic manifold $(N^{2n},\Omega)$ such that $H$ reaches its minimum at $M^{2m}$, a compact, symplectic, Bott-nondegenerate submanifold.

\begin{Theorem}
\labell{Theorem:main} Let $d^2_N H$ satisfy the global resonance conditions as defined below. Then, for sufficiently small $\epsilon > 0$, the number of periodic trajectories of the Hamiltonian vector field $X_H$ on $\{H=\epsilon^2 \}$ is at least $qCL(M)+(n-m)$. Here $q$ is the number of integer independent eigenvalue classes of $d^2_N H$.

\end{Theorem}

\subsection{The Global Resonance Conditions (G.R.C.)}
\label{subsec:grc}

In \cite{ly}, Lyapunov considers the existence of periodic orbits near critical points of a Hamiltonian $H \colon \R^{2n} \ra \R$. For a local minimum he establishes the existence of $n$ periodic orbits on all sufficiently close energy levels provided that the (imaginary) eigenvalues of the linearized Hamiltonian vector field are integer independent. Weinstein achieves the same lower bound in \cite{we:sta} without any resonance assumptions. In replacing the critical point by a symplectic minimum submanifold $M$, we must readopt some of these resonance conditions globally on $M$. A condition similar to that defined below appears in the work of F. Takens \cite{tak} on strong constraints.

 Setting $H|_M=0$, we can assume, for sufficiently small $\epsilon$, that $\{H=\epsilon^2\}$ lies in a neighborhood $W$ of the zero section in the total space $E$ of a normal bundle over $M$. We choose this normal bundle to be $(TM)^{\Omega}$, the symplectic orthogonal complement to $TM$. By Weinstein's Symplectic Neighborhood Theorem \cite{we:lag} we may also assume that that $\Omega$ restricts to the fibres in $(TM)^{\Omega}\cap W$ as a constant linear symplectic form $\Omega^F$. Thus we have a symplectic vector bundle and a fibrewise postitive-definite quadratic form $d^2_N H$. In each fibre $E_x$, there exist coordinates $\{y_i\}_{i=1}^{2p}$ such that $\Omega^F$ is the canonical symplectic form on ${\R}^{2p}$ and $$d^2_N H(x)(y)= \sum_{i=1}^p a_i(x)(y_i^2+y_{i+p}^2).$$ This follows from the fact that one can put a symplectic matrix in standard form and diagonalize a symmetric nondegenerate matrix simultaneously. In fact, the $a_i(x)$ are uniquely determined by $\Omega^F(x)$ and $d^2_N H(x).$  

\begin{Definition}
\labell{Definition:const}
For a fixed $x_0 \in M$ let the $\{a_i\}$ be ordered by magnitude i.e. $$0<a_1 \leq a_2 \ldots \leq a_p.$$ We say that $d^2_N H$ satifies the global resonance conditions if the $a_i$ satisfy the following restrictions. If $a_i(x_0) \neq na_j(x_0)$ for $i > j$ and $ n \in \N$ then $a_i(x) \neq na_j(x)$ for all $x \in M$. Similarly, if $a_i(x_0) = na_j(x_0)$ for $i > j$ and any $ n \in \N$ then $a_i(x)= na_j(x)$ for all $x \in M$. We call such pairs of eigenvalues integer independent and dependent respectively and, ignoring our imposed ordering by magnitude, we denote the classes of integer dependent eigenvalues as $$\{(a_{1}, a_{2},\ldots ,a_{k_1}), (a_{k_{1}+1},\ldots ,a_{k_2}),\ldots, (a_{k_{q-1}+1},\ldots, a_{k_q}=a_p)\}.$$ Under these assumptions the $a_i(x)$ can be treated as smooth real-valued functions of $M$.

\begin{Remark}
\labell{Remark:grc}
The existence of periodic orbits can be established in a variety of situations even if the G.R.C. are not strictly satisfied. For example, we will call a subset of the $\{a_i\}$'s, $(a_{j_1},\ldots,a_{j_r})$, a stable eigenvalue set of order $r$ if the $\{a_{j_i}\}$ are integer dependent along $M$ and no other eigenvalue is an integer multiple of any of the $a_{j_i}$ for any $x \in M$. Under the other hypotheses of Theorem \ref{Theorem:main} the existence of a stable eigenvalue set of order $r$ yields $CL(M)+r$ periodic orbits on sufficiently low energy levels.

Alternatively, the G.R.C. are satisfies automatically if $M$ is of codimension two in $N$. In such cases there is only one eigenvalue, $a_1$, and we establish the existence of at least $Crit(M)$ periodic orbits on low energy levels. Again, $Crit(M)$ denotes the minimum number of critical points of a smooth real-valued function on $M$.   
\end{Remark}

\end{Definition}

\subsection{The Magnetic Problem}

Theorem \ref{Theorem:main} relates to the magnetic problem in the following way. Let $(M,\omega)$ be a compact symplectic manifold and $N=T^* M$ with the symplectic form $\Omega=d\lambda + \pi^* \omega$, where $\lambda$ is the standard Liouville form on $T^* M$ and $\pi:T^*M \ra M$ is the canonoical projection. Then any Riemannian metric $g$ on $M$ will yield a proper kinetic energy Hamiltonian on $T^*M$ for which $M$ is a Bott-nondegenerate minimum. The resulting Hamiltonian vector field on $T^*M$ describes the motion of a particle, constained to $M$, under the influence of the magnetic field given by $\omega$. Finally, $M$ is a symplectic submanifold of $T^*M$ when identified with the zero section and Theorem \ref{Theorem:main} now yields: 

\begin{Theorem}

A charged particle on a manifold $M^{2m}$, under the influence of a ``symplectic'' magnetic field, with a kinetic energy satisfying the G.R.C. has at least $qCL(M)+m$ closed trajectories on sufficiently low energy levels.   

\end{Theorem}

For the magnetic problem the G.R.C. are satisfied automatically in several interesting cases. If $M$ is a surface there is only one $a_i$ (see Remark \ref{Remark:grc}). In this case we get an alternative proof of the following result of Ginzburg in \cite{gi:FA}.

\begin{Theorem}
\labell{Theorem:giFA}

A charged particle on an orientable compact surface $M$, under the influence of a nonvanishing tranversal magnetic field has at least $Crit(M)$ closed trajectories on sufficiently low energy levels.

\end{Theorem}
 
As well, if $g=\omega(\cdot,J\cdot)$  for a compatible almost complex structure $J$ on $(M,\omega)$ we shall see (Remark \ref{Remark:kahler}) that not only are the G.R.C. satisfied but $q=1$ and all the $a_i$ are equal. This includes the case of a Kahlerian metric on a Kahler manifold and under this assumption we get:

\begin{Corollary}
\labell{Theorem:compatible}   
A charged particle on $(M,\omega,J)$, under the influence of the ``symplectic'' magnetic field, has at least $CL(M)+m$ closed trajectories on sufficiently low energy levels.
\end{Corollary}
When the orbits here are nondegenerate the lower bound here may be replaced by $SB(M)$, the sum of Betti numbers of $M$. This is shown in \cite{gi:ker} using the methods of \cite{gi:FA}. 

Finally, as noted by A. Weinstein, if a Lie group acts transitivley on $M$ such that $\omega$ and $g$ are invariant under the action then the G.R.C. are again satisfied. In fact, by the transitivity of the action and the invariance of the vector field we get at least one periodic orbit through every point on $M$.
  
\section{Proof of the Main Theorem}
\labell{sec:proof}

Theorem \ref{Theorem:main} will follow from two lemmae. We first rescale the problem and show that $X_1$, the pushforward of $X_H$ after the rescaling, is close to a fibrewise quasiperiodic vector field $X_0$ on low energy levels. Then with the second lemma we establish the desired number of periodic orbits for $X_1$ and hence $X_H$. 

In looking for periodic orbits near critical points of a  Hamiltonian  $H \colon \R^{2n} \ra \R$, one starts with an equation of the form $$ \dot{z}=J\,H_{zz}(0)\,z+R(z)$$ where $\small{R(0)=0=\frac{\partial R}{\partial z}(0)}$ and 
$\scriptsize{  J = \left[ \begin{array}{cc}
          0     & -Id_n\\
          Id_n  & 0  
               \end{array}
         \right]}. $
This equation is then rescaled using $z \ra \epsilon z$, i.e. $$ \dot{z}= J\,H_{zz}(0)\,z+O(\epsilon),$$ and one looks for periodic orbits near those of the limiting linear equations,(\cite{ly}, \cite{we:sta}, \cite{Moser-1976}, \cite{fara}). In dealing with our critical submanifold $M$ we will rescale globally in the normal directions of the tubular neighborhood of $M$ described in \ref{subsec:grc} and show that we still get a well-defined and useful limiting vector field.
 
We start with the Hamiltonian dynamical system defined by $$i_{X_H}\Omega=dH.$$
Letting $\Phi\colon E \ra E$ be the global fibrewise dilation by a factor of $\epsilon$, we set $$X_1=\epsilon^{-2}{\Phi^{-1}}_* X_H,$$ $$\tilde{\Omega}= \Phi^*\Omega,$$ and $$\tilde{\Omega}^F=\Phi^*\Omega^F.$$ The Hamiltonian dynamical system now looks like:
\begin{equation}
\labell{eq:dyn}
i_{X_1}\tilde{\Omega}=d(\epsilon^{-2}\Phi^* H).
\end{equation} 

\begin{Lemma}
\labell{Lemma:lemma1}

As $\epsilon \ra 0$, $X_1$ approaches a fibrewise, quasiperiodic vector field $X_0$. In particular, in each fibre $E_x$, $X_0$ is the Hamiltonian vector field of $d^2_N H|_x$ with respect to $\tilde{\Omega}^F$. 

\end{Lemma}

\begin{proof}

We restrict ourselves to a neighborhood of $W$ on which our normal bundle is trivial and there we split things into fibre (vertical) and base (horizontal) components, i.e. let $$X_1=X^B+X^F$$ and $$\tilde{\Omega}= \omega +\eta + \xi$$ where $\omega$ has only horizontal terms, $\xi$ has only fibre terms and $\eta$ has only mixed terms. The forms $\omega$ and $\xi$ have orders of $\epsilon^{-2}$ and 1 respectively. Furthermore, by our choice of $(TM)^{\Omega}$ as the normal bundle, we have $\xi =\tilde{\Omega}^F$ and $\eta |_{TM}=0$ so that $\eta$ is of (at most) order 1 in $W$. Locally, Hamilton's equations \eqref{eq:dyn} are now of the form
\begin{eqnarray*}
i_{X^B} \omega + i_{X^F}\eta           & = & d^B(\epsilon^{-2}\Phi^*H)\\
i_{X^F}\tilde{\Omega}^F  + i_{X^B}\eta & = & d^F(\epsilon^{-2}\Phi^*H).
\end{eqnarray*}
Comparing relative orders, the first equation implies that $X^B$ is of order $\epsilon^2$ with respect to $X^F$. The second equation then implies that $X^F$ is of order 1. We also know that, as $\epsilon \ra 0$,  $\epsilon^{-2}\Phi_*H$ approaches $d^2_N H$. Hence, in this limit,  $X_1$ approaches the fibrewise vector field $X_0$ given by
\begin{equation}
\labell{eq:x0} 
i_{X_0}\tilde{\Omega}^F=d^F( d^2_N H).
\end{equation}
Indeed, this equation defines the limiting vector field $X_0$ globally and we note that the convergence is $C^k$ for any $k$.

To show that the flow of $X_0$ is fibrewise quasiperiodic, we recall that in each fibre $E_x$ there are coordinates $\{y_i(x)\}$ such that $$\tilde{\Omega}^F=\sum_{i=1}^{k_q} dy_i \wedge d y_{i+k_q}$$ and $$d^F(d^2_N H)= 2 \sum_{i=1}^{k_q} a_i(x)(y_i dy_i+ y_{i+k_q} dy_{i+k_q}).$$
So, the flow of $X_0$ in each fibre is given by the sets of O.D.E's
\begin{eqnarray*}
\dot{y}_i        & = & -2a_i y_{i+k_q}\\
\dot{y}_{i+k_q}  & = & 2a_i y_i                
\end{eqnarray*}
 for $i \in \{1,\ldots,k_q=(n-m)\}.$ This flow is clearly quasiperiodic and the proof of the lemma is done.
\end{proof}

\begin{Remark}
\labell{Remark:kahler}As an aside, we show that all the $\{a_i\}$ are equal for the magnetic problem when the metric $g$ comes from an almost complex structure $J$ compatible with $\omega$. In the framework of the proof of Lemma \ref{Lemma:lemma1} we have $N=T^*M$, $\Omega= \pi^*{\omega}+d\lambda$ and $H(q,p)=\haf \omega(q)(p,J(p))$, where $\{q,p\}$ are local ``$pdq$'' coordinates of $\lambda$. It is possible to choose the base coordinates $\{q_i\}$ such that $\omega$ is locally canonical. Thinking of $\Omega$ in these coordinates as a skew symmetric matrix we then have   
\[          \Omega = \left[ \begin{array}{cc}
                               \omega & -Id_m\\
                               Id_m   & 0  
                          \end{array}
                     \right].
                                \]
Using the normal bundle $(TM)^{\Omega} \subset T(T^*M)$, it is staight forward to check that ${\Omega}^F=-\omega$, as matrices. After rescaling, we see that as $\epsilon \ra 0$ the fibres of $(TM)^{\Omega}$ approach the vertical fibres of $T(T^*M)$ and there we get the limiting vector field $X_0$ given by $$-\omega X_0(q,p)=\omega J(p).$$ The nondegeneracy of $\omega$ implies that $X_0(q,p)=-J(p)$ and choosing the $\{p_i\}$ such that $J$ is canonical in $T^*_q M$ we see that all the $a_i=1$. 
\end{Remark}

Next we show that the flow of $X_0$, $\phi^t$, has $q$ nondegenerate periodic submanifolds in $\{ d^2_N H=1\}$. Recall that the G.R.C. yield $q$ classes of integer dependent eigenvalues of $d^2_N H$ with respect to $\tilde{\Omega}^F$: 
$$\{(a_{1}, a_{2},\ldots ,a_{k_1}), (a_{k_{1}+1},\ldots ,a_{k_2}),\ldots, (a_{k_{p-1}+1},\ldots, a_{k_p})\}.$$ By the proof of the previous lemma we know that the flow of $X_0$ resticted to the eigenspace of any of these classes is periodic. The level set $\{d^2_N H=1\}$ is a sphere bundle over $M$. Consider the subbundles, $\{\Sigma^{i}\}_{i=1}^q$, which correspond to the $q$ classes of integer dependent eigenvalues of $d^2_N H$ as follows.  Let each $\Sigma^{i}$ be the $S^{2(k_i-k_{i-1})-1}$ bundle over $M$ whose fibre over $x \in M$ is given by the intersection of the fibre of $\{d^2_N H=1\}$ with the eigenspace of the class $(a_{k_{i-1}},\ldots,a_{k_i})$. By the integer independence of the classes, the $\Sigma^{i}$ are nondegenerate as periodic submanifolds of the flow of $X_0$ on $\{d^2_N H=1\}$.

We now look for closed orbits of $X_1$ near the $\Sigma^i$. 

\begin{Lemma}
\labell{Lemma:lemma2}
 
The Hamiltonian vector field $X_1$ has at least $Crit(\Sigma^{i}/S^1)$ closed orbits near each $\Sigma^{i}$.

\end{Lemma}

\begin{proof}

The proof will come in two steps, both dependent on Bottkol's generalizations of Moser's work. Fix $\Sigma^{i}=\Sigma$. First we will use a proposition of Bottkol's to find a $\phi^t$-invariant vector field  $V$ on $\Sigma$ whose zeroes correspond to closed orbits of $X_1$ near $\Sigma$. Then we will extend the method of proof of Bottkol's second result in \cite{Bottkol} to find a $\phi^t$-invariant function $S$ on $\Sigma$ whose critical points correspond  directly to zeroes of $V$. Combining these results we get the following as desired:
\begin{eqnarray*}
\{ \#\;of\;closed\;orbits\;of\;X_1\;near\;\Sigma\} & \geq & \{\#\;of\; zeroes\; of\; V\; on\; \Sigma\}\\
                                                   & \geq & Crit(\Sigma / S^1)
\end{eqnarray*}
\subsection{Step 1}

In \cite{Bottkol}, the following local decomposition is proved.

\begin{Proposition} 
\labell{Proposition:bottkol1}
Let $Y_0$ be a $C^3$ vector field on a manifold $P$ that generates a flow $\phi^t$ with a compact nondegenerate periodic submanifold $\Sigma$ of period 1. Suppose that $Y_1$ is a $C^2$ vector field $C^2$-close to $Y_0$ in some neighborhood of $\Sigma$. Then there exists a ${\phi}^t$-invariant function $\lambda \colon \Sigma \ra {\R}$ (close to 1),  and unique $C^1$-small sections U and V, of ${TP|}_{\Sigma}$ and $T\Sigma$ respectively, such that for $$ u=exp(U(\cdot)): \Sigma \ra P$$ and $$P(u)\colon T_z P \ra T_{u(z)} P,$$ defined by $$P(u)V={ \frac{d}{dh}\ \biggl\vert}_{h=0} exp(U+hV),$$ we have the following decomposition:

$$ \lambda Y_1 (u)= duY_0-P(u)V. $$

The following conditions also hold:

\begin{enumerate}

\item 
\emph{} $[V,Y_0]=0$.

\item 
\emph{} $\int_0^1 d \phi^{-t} U_{\Sigma}(\phi^t (z))\, dt=0 \;\;\; \forall z \in  \Sigma$.

\item 
\emph{}$V \perp Y_0$ with respect to a $\phi^t$-invariant metric on $\Sigma$, $<\;,\;>$.

\end{enumerate}

\end{Proposition}

We now outline a simple proof of this proposition for the case $\Sigma=P$. Consider the following subsets of $\chi^1(P)$, the Banach space of $C^1$ vector fields on $P$. Let $$\tilde{\chi}^1(P)=\{U \in \chi^1 (P)\; \bigl\vert\; [Y_0,U] \in \chi^1 (P) \; \; and\; \; \int_0^1 d\phi^{-t} U(\phi^t (z)) \,dt=0 \;\; \forall z \in P \}$$ and $$\chi^1_{inv}=\{V \in \chi^1(P) \; \bigl\vert \; [Y_0,V]=0\}.$$ Clearly $\chi^1_{inv}$ is a closed subspace of $\chi^1(P)$, and $\tilde{\chi}^1(P)$ is also a Banach space given the norm $$ \|U\|_{~}=\|U\|_1+\|\,[Y_0,U]\,\|_1.$$ Let the map $\Pi \colon \tilde{\chi}^1(P) \times \chi^1_{inv} \ra \chi^1(P)$ be defined by $$\Pi(U,V)=duY_0-P(u)V,$$ where $u$ and P(u) are as above. Note that $\Pi(0,0)=Y_0$. One can show quite easily that $d\Pi_{(0,0)}$ is an isomorphism. So, by the inverse function theorem for Banach spaces, for all $Y_1$ sufficienly $C^1$-close to $Y_0$ we get unique $U$ and $V$ such that $$Y_1=duY_0-P(u)V.$$ Given this decomposition, the existence of the $\phi^t$-invariant function $\lambda$, that ensures $V \perp Y_0$, can be established using the implicit function theorem  as it is in \cite{Bottkol}.

\begin{Remark} With a little extra effort this proof yields a $C^1$ version of a result of H. Seifert that establishes the existence of a closed orbit for vector fields $C^0$-close to the Hopf field on $S^3$, \cite{sei}. It is a simplification of the proof by Moser, \cite{Moser-1976}, which is possible only because of its stronger closeness hypothesis, (see \cite{gi:per}).
\end{Remark}
Proposition \ref{Proposition:bottkol1} yields a decomposition of $Y_1$ near $\Sigma$, i.e. there is a unique embedding $u$ of $\Sigma \ra P$, and a $\phi^t$-invariant vector field $V$ on $\Sigma$ such that on $u(\Sigma)$ we have $$ \lambda Y_1(u(z))=duY_0-P(u)V. $$ This is a useful splitting because if $V(z)=0$ then $V(\phi^t (z))=0$ for all $t$, since $[V,Y_0]=0$. Thus, we get $$ \lambda Y_1(u(\phi^t (z)))=duY_0(\phi^t(z))$$ which means that $u(\phi^t (z))$ is a closed orbit of $ \lambda Y_1$ on $u(\Sigma).$ Since $\lambda$ is close to 1, $u(\phi^t (z))$ is also a closed orbit of $Y_1$ and we get a direct correspondence between zeroes of $V$ and periodic orbits of $Y_1$ on $u(\Sigma)$.

In looking to apply this decomposition to our situation we note, for small $\epsilon >0$, that $\epsilon^{-2}\Phi_*H$ is $C^k$-close to $d^2_N H$ for any $k$. Consequently there exists a diffeomorphism $\beta \colon \{d^2_N H=1\} \ra \{\epsilon^{-2}\Phi_*H=1\}$ such that $d\beta$ is close to the identity. Letting $\{d^2_N H=1\}=P$, $\Sigma =\Sigma$, $X_0=Y_0$, and $Y_1=\beta^* X_1$, Proposition \ref{Proposition:bottkol1} gives unique $U$ and $V$ such that
\begin{equation}
\labell{eq:dec}
\lambda \beta^* X_1(u(z))=duX_0-P(u)V. 
\end{equation}
Now if $V(z)=0$ for $z \in \Sigma$ we get $$ \lambda \beta^* X_1 (u(\phi^t(z)))=duX_0(\phi^t(z)).$$ This implies that $\beta \circ u \circ \phi^t (z)$ is a closed orbit of $X_1$ on $u(\Sigma) \subset \{\epsilon^{-2}\Phi_*H=1\}$ and we get the first desired correspondence between zeroes of a $\phi^t$-invariant vector field $V$ on $\Sigma$ and closed orbits of $X_1$ near $\Sigma$.

\subsection{Step 2}

Expanding on the work of Bottkol we constuct a $\phi^t$-invariant function $S \colon \Sigma \ra \R$ whose critical points correspond to zeroes of $V$. Define the following 1-form on $\Sigma$ 
\begin{equation}
\labell{eq:S}
dS=\int_0^1 \gamma_t^*(i_{\dot{\gamma}_t} \tilde{\Omega})\,dt ,
\end{equation}
where $\gamma_t(z)=\beta \circ u \circ \phi^t(z)$.
\begin{Claim}
As suggested by the notation, $dS$ is exact.
\end{Claim}

\begin{proof}
From the definition of $\gamma_t$ we have 
\begin{equation*}
dS = \int_0^1 {\phi^t}^*(i_{ X_0} u^*\beta^*\tilde{\Omega})\,dt. 
\end{equation*}
The map $\beta \circ u$ is homotopic to the identity map on $\Sigma \subset W \subset E$. By the homotopy invariance of the de Rham cohomology we get $u^*\beta^*\tilde{\Omega}=\tilde{\Omega}+d\lambda$ and 
\begin{equation}
\labell{eq:exa}
dS= \int_0^1 {\phi^t}^*(i_{ X_0}\tilde{\Omega})\,dt +  \int_0^1 {\phi^t}^*(i_{ X_0} d\lambda)\,dt.
\end{equation}
The second term of \eqref{eq:exa} is exact since
\begin{eqnarray*}
\int_0^1 {\phi^t}^*(i_{ X_0} d\lambda)\,dt & = & \int_0^1 {\phi^t}^*(L_{ X_0}\lambda)\,dt-\int_0^1 {\phi^t}^*(di_{X_0} \lambda)\,dt\\
                                           & = & d\int_0^1 {\phi^t}^*(i_{X_0}\lambda)\,dt.
\end{eqnarray*}
 Here the term involving $L_{X_0}$ vanishes after integrating over the closed curve.

We still must show that the first term in \eqref{eq:exa} is exact. Overall, we have the following situation $$E \stackrel{I_{\Sigma}}{\hookleftarrow} \Sigma \stackrel{p_1}{\longrightarrow} \Sigma/S^1 \stackrel{p_2}{\longrightarrow} M$$ and $$\pi \colon E \ra M,$$ where $\pi$, $p_1$ and $p_2$ are the obvious projections and $I_{\Sigma}$ is inclusion. Again by homotopy invariance, we have $[\tilde{\Omega}] = [\pi^* \sigma]$ for some 2-form $\sigma$ on $M$. We also know that $$\pi \circ I_{\Sigma} = p_2 \circ p_1.$$ Hence
\begin{eqnarray*}
[I_{\Sigma}^*\tilde{\Omega}]  & = & [I_{\Sigma}^* \circ \pi^* (\sigma)]\\
                              & = & [p_1^* \circ p_2^* (\sigma)]
\end{eqnarray*}
and we see that $[I_{\Sigma}^* \tilde{\Omega}] = [p_1^*\mu]$ for some 2-form $\mu$ on $\Sigma/S^1$. This means we can write $\tilde{\Omega}= \Omega_0 + d\nu$ where $\Omega_0$ satisfies $i_{X_0} \Omega_0=0$. Consequently
\begin{eqnarray*}
\int_0^1 {\phi^t}^*(i_{ X_0}\tilde{\Omega})\,dt  & = & \int_0^1{\phi^t}^*(i_{ X_0}\Omega_0)\,dt+\int_0^1 {\phi^t}^*(i_{ X_0}d\nu)\,dt\\
                                                 & = & \int_0^1{\phi^t}^*(L_{X_0}\nu)\,dt - \int_0^1 {\phi^t}^*(d i_{ X_0}\nu)\,dt\\
                                                 & = & - d\int_0^1 {\phi^t}^*(i_{ X_0}\nu)\,dt.
\end{eqnarray*}
This completes the proof of the claim. 
\end{proof}
\begin{Claim} $dS$ is a $\phi^t$-invariant form.
\end{Claim}

\begin{proof}\begin{eqnarray*}
dS(X_0) & = &  i_{X_0}\int_0^1 \gamma_t^* \tilde{\Omega}(\dot{\gamma}_t , \cdot )\,dt \\
        & = & \int_0^1 \tilde{\Omega}(\dot{\gamma}_t,{\gamma_t}_*X_0)\,dt 
\end{eqnarray*}
Using the definition of $\gamma_t$ we have $$ \dot{\gamma}_t(z)=d\beta \circ du X_0(\phi^t(z))$$ and $${\gamma_t}_*=d\beta \circ du \circ d\phi^t.$$ Hence
\begin{eqnarray*}
dS(X_0)(z) & = & \int_0^1 \tilde{\Omega}(d\beta \circ du X_0(\phi^t(z)) ,d\beta \circ du \circ d\phi^t X_0(z) )\,dt \\
           & = & 0
\end{eqnarray*}
\end{proof}
The previous two claims establish the existence of a $\phi^t$-invariant function $S$ whose differential is given by \eqref{eq:S}.

\begin{Remark}
\labell{Remark:ind}
 It is worthwhile to note here the following two facts. First, our construction of $S$ does not depend on the relation of the embeddings $u$ and $\beta$ to the vector fields $X_0$ and $X_1$. In fact, any $u$ and $\beta$ will yield a $\phi^t$-invariant function $S$ in the same manner. This freedom of choice will be used in what follows. Also, just as we used the exponential map to obtain the diffeomorphism $u$ from the section $U$, we can consider $\beta$ as being defined by a small smooth section $B$ of ${TE|}_{\{d^2_N H=1\}}$, i.e. $\beta (z)=exp(B(z)).$ 
\end{Remark}

We still need to show that the critical points of $S$ correspond to zeroes of $V$. Using our expression for $\dot{\gamma}_t(z)$ and the decomposition \eqref{eq:dec} we get $$ \dot{\gamma}_t(z)= d\beta P(u)d\phi^t V(z) - \lambda (\phi^t (z) )X_1 (\beta \circ u (\phi^t(z))).$$ Substituting this into the definition of $dS$ yields
\begin{equation*}
dS=\int_0^1 \gamma_t^* \tilde{\Omega} (d\beta P(u)d\phi^tV(z),\cdot)\,dt-\int_0^1 \gamma_t^* \tilde{\Omega} (\lambda (\phi^t (z) )X_1 (\beta \circ u (\phi^t(z))),\cdot)\,dt.
\end{equation*}
The term involving $X_1$ vanishes since $i_{X_1}\tilde{\Omega}=0$ on $\{\epsilon^{-2}\Phi_*H=1\}$ and we are left with
\begin{equation}
\labell{eq:rel}
dS=\int_0^1 \gamma_t^* \tilde{\Omega} (d\beta P(u)d\phi^tV(z),\cdot)\,dt.
\end{equation}
Since $dS(X_0)=0$ we know that $dS$ is a 1-form on ${X_0}^{\perp} \subset T\Sigma$ where $\perp$ is taken with respect to a fixed $\phi^t$-invariant metric $<\;,\;>$. Fixing a $z \in \Sigma$ we can consider equation \eqref{eq:rel} as defining a family of linear maps  $$K(U,B) \colon ({X_0}^{\perp})_z \ra  ({X_0}^{\perp})_z^*$$ parameterized by the embeddings $u$ and $\beta$ and hence the sections $U$ and $B$, (see Remark \ref{Remark:ind}). Here $ ({X_0}^{\perp})_z^*$ is the dual space to ${X_0}^{\perp}$ at $z$ and $K(U,B)(V)=dS$, as defined in \eqref{eq:rel}.

We know that $U$ and $B$ are small in their $C^1$-norms. If we can show that $K(0,0)$ is nondegenerate then for our small sections, $U$ and $B$, $K(U,B)$ will also be nondegenerate. This would imply that $$ dS=0 \iff V=0 $$ and we would have the desired correspondence between the zeroes of $V$ and the critical points of the $\phi^t$-invariant function $S$ on $\Sigma$.
  
\begin{Claim}

$K(0,0)$ is nondegenerate.

\end{Claim}

\begin{proof}

With $U=0$ and $B=0$ we have $\gamma_t(z)=\phi^t(z)$ and $K(0,0) \colon X_0 \ra (X_0)^{\perp}$ is given by $$dS=\int_0^1 ({\phi^t}^* \tilde{\Omega})(V,\cdot)dt.$$

The 2-form $\tilde{\Omega}$ is almost conserved by ${\phi^t}^*$ since $\phi^t$ is close to the flow of $X_1$. So we have $$ {\phi^t}^*\tilde{\Omega} = \tilde{\Omega} + \underbrace{E(t)}_{order \; \epsilon}$$ and $$dS=\tilde{\Omega}(V,\cdot)+\underbrace{\int_0^1E(t)(V,\cdot)\,dt}_{order \; \epsilon}.$$ All that is required now is to show that ${\tilde{\Omega}|}_{X_0^{\perp}}$ is nondegenerate. 
The kernel of $\tilde{\Omega}$ restricted to $T\Sigma$ is $X_{1,\Sigma},$ the projection of $X_1$ onto $T\Sigma$. But $X_0$ is $C^k$-close to $X_1$ and hence to $X_{1,\Sigma}$. The nondegeneracy of $\tilde{\Omega}$ on $(X_{1,\Sigma})^{\perp}$ implies that it is also nondegenerate on $(X_0)^{\perp}$ and both the claim and Step 2 follow. 
   
\end{proof}
 
For each $\Sigma^i$ we have constructed a $\phi^t$-invariant function $S_i$ whose critical submanifolds correspond to distinct periodic trajectories of $X_1$ near $\Sigma^i$. Now, $\phi^t$ generates a circle action on each $\Sigma^i$ without fixed points. This yields $$\sum_{i=1}^q Crit(\Sigma^i/S^1)$$ as an immediate lower bound for the number of closed orbits on low energy levels. If all the $\{a_j\}_{j=k_{i-1}}^{k_i}$ are equal then the circle action generated by $\phi^t$ on $\Sigma^i$ is free. In particular, $\Sigma^i/S^1$ is a $\C P^{(k_i-k_{i-1}-1)}$ bundle over $M$ obtained from a vector bundle. For such bundles we have that
\begin{equation}
\labell{eq:cohom}
H^*(\Sigma^i/S^1;\R)=H^*(M;\R) \otimes H^*(\C P^{(k_i-k_{i-1}-1)};\R)
\end{equation}
as modules over $H^*(M;R)$, \cite{hus}.
Using Ljusternik-Schnirelman theory we then get \begin{eqnarray*}
Crit(\Sigma^i/S^1)   &\geq& Cat(\Sigma^i/S^1)\\
                     &\geq& CL(\Sigma^i/S^1)+1\\
                     & =  & CL(M)+k_i-k_{i-1}.
\end{eqnarray*}
In fact, even if the ${a_i}$ aren't equal, and $\Sigma^i/S^1$ is a twisted projective space bundle, it was shown by Weinstein in \cite{we:v} that
\begin{equation*}
\{\#\;of\;crit.\;circles\;of\;S_i\;on\;\Sigma^i \} \geq Cat(\Sigma^i/S^1).
\end{equation*}
Relation \eqref{eq:cohom} still holds in this case and so we always have 
\begin{eqnarray*}
\{ \#\;of\;periodic\;orbits\;of\;X_1\;near\;\Sigma^i\}  &\geq& CL(M)+k_i - k_{i-1}.
\end{eqnarray*}
Summing over the $\Sigma^i$ and using the fact that $k_0=0$ and $k_q=(n-m)$ we get at least $qCL(M)+(n-m)$ distinct closed orbits of $X_1$ and hence $X_H$. The proof of Theorem \ref{Theorem:main} is now complete.

\end{proof}

As a point of further interest we note that the proof of Lemma \ref{Lemma:lemma2} encompasses a proof of the following more general result. Let $Q$ be an odd dimensional manifold and $X_0$ a nonvanishing $C^3$ vector field on $Q$ whose flow has a compact nondegenerate periodic submanifold $\Sigma.$

\begin{Theorem}
Let $\Omega$ be a closed, maximally nondegenerate 2-form on $Q$ whose kernel is $C^2$-close to the line bundle spanned by $X_0$ near $\Sigma$. If $ \Omega=p_1^*\sigma$ for the canonical projection $p_1 \colon \Sigma \ra \Sigma/S^1$ and some $\sigma \in H^2(\Sigma/S^1)$, then $\Omega$ has at least $Crit(\Sigma/S^1)$ closed characteristics near $\Sigma$. 
\end{Theorem}


\begin{thebibliography}{ABCD}

\bibitem[Ar1]{Arn1}
Arnold, V.I., First steps of symplectic topology, {\em Russian Math. Surveys}, {\bf 41} (6) (1986), 1-21. 

\bibitem[Ar2]{Arn2}
Arnold, V.I., On some problems in symplectic topology, in {\em Topology and Geometry-Rochlin Seminar}, O.Ya. Viro (Editor), Lect. Notes in Math., vol. 1346, Springer, 1988. 

\bibitem[BT]{bata}
Bahri, A., Taimanov, I.A., Periodic orbits in magnetic fields and Ricci curvature of Lagrangian systems, {\em Trans.  A.M.S.}, {\bf 350} (7) (1998), 2697-2717.

\bibitem[Bot]{Bottkol}
Bottkol, M., Bifurcation of periodic orbits on manifolds and
Hamiltonian systems. \emph{J. Diff. Eq.}, {\bf 37} (1980), 
12--22. 

\bibitem[FR]{fara}
Fadell, E. R., Rabinowitz, P. H., Generalized cohomological index theories for Lie group actions and an application to bifurcation problems for Hamiltonian systems, {\em Inv. Math.}, {\bf 45} (1978), 48-67.

\bibitem[FHV]{fhv}
Floer, A., Hofer, H., Viterbo, C., The Weinstein conjecture on $P \times \C^l$, {\em Math. Z.}, {\bf203} (1990), 469-482.

\bibitem[Gi1]{gi:FA}
Ginzburg, V. L., 
New generalizations of Poincar\'{e}'s geometric theorem,
\emph{Funct. Anal. Appl.}, {\bf 21} (2) (1987), 100--106.

\bibitem[Gi2]{gi:Cambr}
Ginzburg, V. L., On closed trajectories of a charge in a magnetic
field. An application of symplectic geometry, in {\em Contact
and Symplectic Geometry}, C.B. Thomas (Editor), Publications of
the Newton Institute, Cambridge University Press, Cambridge, 1996,
p. 131--148.

\bibitem[Gi3]{gi:MathZ}
Ginzburg, V. L., On the existence and non-existence of closed 
trajectories for some Hamiltonian flows, {\em Math. Z.},
{\bf 223} (1996), 397--409.

\bibitem[Gi4]{gi:Glik}
Ginzburg, V.L., Accesible points and closed trajectories of mechanical systems, in {\em Global Analysis in Mathematical Physics} by Yu. Gliklihk, Springer-Verlag, 1998.

\bibitem[Gi5]{gi:per}
Hamiltonian dynamical systems without periodic orbits, Preprint, 1998.

\bibitem[GK]{gi:ker}
Ginzburg, V. L., Kerman, E., Periodic orbits in magnetic fields in dimensions greater than two, Preprint, 1999.

\bibitem[He]{hed}
Hedlund, G.A., Fuschian groups and transitive horocycles, {\em Duke Math. J.}, {\bf 2} (1936), 530-542.

\bibitem[Ho]{hof1}
Hofer, H., On the topological properties of a symplectic map, {\em Proc. Royal Soc. Edinburgh}, {\bf 115A} (1990), 25-38.

\bibitem[HZ]{hz}
Hofer, H., Zehnder, E., {\em Symplectic invariants and Hamiltonian dynamics}, Birkh\"{a}user, Basel, 1994.

\bibitem[Hu]{hus}
Husemoller, D., {\em Fibre Bundles}, Springer-Verlag, New York, 1966.

\bibitem[Ji]{jiang}
Jiang, M.-Y., Hofer-Zehnder symplectic capacity for 2-dimensional manifolds, {\em Proc. Roy. Soc. Edinburgh}, {\bf123A} (1993), 945-950.

\bibitem[Ly]{ly}
Lyapunov, A., Probl\`{e}me g\'{e}n\'{e}rale de la stabilite du mouvement, {\em Ann. Fac. Sci Toulouse}, {\bf2} (1907), 203-474.

\bibitem[Ma]{ma}
Mather, J.N., Private communication to V. Ginzburg, November 1995.

\bibitem[Mo1]{moser:kep}
Moser, J., Regularization of Kepler's problem and the averaging method on a manifold, {\em Comm. Pure Appl. Math.}, {\bf27} (1970), 609-636.

\bibitem[Mo2]{Moser-1976}
Moser, J., Periodic orbits near an equilibrium and a Theorem
by Alan Weinstein, \emph{Comm. Pure and Appl. Math.}, {\bf 29}
(1976), 727--747.

\bibitem[Po]{pol}
Polterovich, L., Geometry on the group of Hamiltonian Diffeomorphisms, {\em Procedings of the I.C.M.}, {\rm Vol.II} (Berlin 1998), 1998, 401-410.

\bibitem[Se]{sei}
Seifert, H., Closed integral curves in 3-space and isotopic two-dimensional deformations, {\em Proc. Amer. Math. Soc.}, {\bf1} (1950), 287-302.

\bibitem[Ta]{tak}
Takens, F., Motion under the influence of a strong constraining force, {\em Springer lecture notes in math.}, Vol. 819, Springer, Berlin, 1980, 425-445.

\bibitem[We1]{we:persub}
Weinstein, A., Perturbation of periodic manifolds of Hamiltonian systems, {\em Bulletin Amer. Math. Soc.}, {\bf77} (1971), 814-818.

\bibitem[We2]{we:lag}
Weinstein, A., Symplectomorphisms and their Lagrangian Submanifolds, {Advances in Mathematics}, {\bf6} (19710, 329-346.

\bibitem[We3]{we:sta}
Weinstein, A., Normal modes for non-linear Hamiltonian systems,
\emph{Invent. Math.}, {\bf 20} (1973), 377--410.

\bibitem[We4]{we:v}
Weinstein, A., Symplectic V-manifolds, periodic orbits of Hamiltonian systems and the volume of certain Riemannian manifolds,
{\em Comm. Pure Appl. Math.}, {\bf 30} (1977), 265-271.
\end{thebibliography}
\end{document}